\documentclass[12pt,a4paper]{article}
\usepackage{amsmath,amssymb,amsfonts}
\def\Bbb{\mathbb}

\title{\bf The fractional parts of Dedekind sums}

\author{Kurt Girstmair}
\date{}

\makeatletter
\let\@@maketitle=\maketitle
\def\maketitle{\def\thispagestyle##1{\relax}\@@maketitle}
\makeatother
\newtheorem{theorem}{Theorem}

\def\BE{\begin{equation}}
\def\EE{\end{equation}}
\def\BD{\begin{displaymath}}
\def\ED{\end{displaymath}}
\def\BA{\begin{array}}
\def\EA{\end{array}}
\def\BEA{\begin{eqnarray*}}
\def\EEA{\end{eqnarray*}}
\def\BI{\bibitem}

\def\Z{\Bbb Z}

\def\R{\Bbb R}

\def\phi{\varphi}

\def\CMOD#1#2#3{#1 \equiv #2 \: \mbox{mod}\: #3}

\def\MB{\mbox}
\def\LD{\ldots}

\def\NDIV{\, \nmid \,}

\def\BQ{``}
\def\EQ{'' }

\def\DED{Dedekind }

\begin{document}
\maketitle

\begin{abstract}
\noindent
We show that each rational number $r$, $0\le r<1$, occurs as the fractional part of a \DED sum $S(m,n)$.
\end{abstract}

\section*{1. Introduction and result}

Let $n$ be a positive integer and $m\in \Z$, $(m,n)=1$. The classical \DED sum $s(m,n)$ is defined by
\BD
   s(m,n)=\sum_{k=1}^n ((k/n))((mk/n))
\ED
where $((\LD))$ is the \BQ sawtooth function\EQ defined by
\BD
  ((t))=\left\{\begin{array}{ll}
                 t-\lfloor t\rfloor-1/2 & \MB{ if } t\in\R\smallsetminus \Z; \\
                 0 & \MB{ if } t\in \Z
               \end{array}\right.
\ED
(see, for instance, \cite[p. 1]{RaGr}). In the present setting it is more
natural to work with
\BD
 S(m,n)=12s(m,n)
\ED instead.
Observe that $S(m+n,n)=S(m,n)$, so one often considers
only arguments $m$ in the range $0\le m\le n-1$.

\DED sums have quite a number of interesting applications in
analytic number theory (modular forms), algebraic number theory (class numbers),
lattice point problems and algebraic geometry (for instance \cite{{Ap}, {Me}, {RaGr}, {Ur}}).
Moreover, the distribution of these sums has been the subject of study of several authors
(see \cite{{Br}, {Gi}, {Va}, {Zh}}).

The values of \DED sums are rational numbers. It follows from a result in \cite{Hi} that the set
\BD
  \{S(m,n): n>0, 1\le m\le n, (m,n)=1\}
\ED
is dense in the set $\R$ of real numbers. In particular, \DED sums come arbitrarily close to a given rational number.
It seems to be unknown, however, which rational numbers actually occur as the values of  \DED sums $S(m,n)$.
In this note we show that the fractional parts of \DED sums $S(m,n)$ take each possible value.

\begin{theorem} 
\label{t1}

Let the integers $n$ and $q$, $0\le q\le n-1,$ $(q,n)=1$, be given. Then there are integers $m$, $n'$,
$0\le m\le n'-1$, $(m,n')=1$, such that
\BD
  S(m,n')\in \frac qn+\Z.
\ED

\end{theorem} 

Our basic tool is the Barkan-Hickerson-Knuth formula (see, for instance \cite{Hi}). It is a consequence of this formula
that
\BE
\label{1}
  S(m,n)\in \frac{m+m^*}n+\Z,
\EE
where $m^*\in\Z$ is such that $\CMOD{mm^*}1n$ (see, for instance \cite[formula (4)]{Gi2}).
Hence it suffices, for our purpose, to study the behaviour of
$(m+m^*)/n$. Note, however, that, in general,  the fractional part of this number does not run through all possible
values $q/n$, $0\le q\le n-1$, $(q,n)=1$, if $m$ runs through all values $0\le m\le n-1$, $(m,n)=1$. For instance,
it is not hard to see that there are at most $(n+1)/2$ values of fractional parts of $(m+m^*)/n$, if $n\ge 3$
is a prime number, but $n-1$ possible numbers $q/n$. Our proof is based on a suitable extension of the denominator $n$.

\section*{Proof of Theorem \ref{t1}}

We may assume $q\ne 0$ and $n\ge 2$, since $S(m,n')=0$ if $\CMOD{m^2}{-1}{n'}$ (see \cite[p. 28]{RaGr}).

{\em Case 1.} Let $n$ be odd. Let $p$ be a prime number such that
\BE
\label{2}
   \CMOD 2{qp}n.
\EE
In particular, $p\NDIV n$.
In addition, let $\CMOD p14$. By Dirichlets theorem of primes in arithmetic progressions,
such a prime $p$ exists (see, for instance, \cite[p. 103]{Se}).
Since $\CMOD p14$, the number $-1$ is a square mod $p$, hence there is an integer $m$ such that
\BE
\label{3}
  \CMOD{m^2+1}{0}p.
\EE
Since $(n,p)=1$, this number $m$ can be chosen such that
\BD
  \CMOD m1n.
\ED
By (\ref{2}),
\BD
  m^2+1\equiv 1+1\equiv qp \MB{ mod } n.
\ED
Therefore, we may write $m^2+1=qp+kn$, with $k\in \Z$. In view of (\ref{3})
and the fact that $p$ does not divide $n$, it must divide $k$. Thus $m^2+1=qp+k_1np$, $k_1\in \Z$.
Since $(m,np)=1$, there is an integer $m^*$ such that $\CMOD{mm^*}1{np}$. Accordingly, we write
$mm^*=1+k_2np$, $k_2\in\Z$. This gives, on the one hand,
\BE
\label{4}
  m^*(m^2+1)= m+m^*+k_2mnp,
\EE
on the other hand,
\BE
\label{5}
  m^*(m^2+1)=m^*(qp+k_1np)=m^*qp+m^*k_1np.
\EE
Since $\CMOD m1n$, we have $\CMOD{m^*}1n$, and so $m^*=1+k_3n$, $k_3\in \Z$.
Accordingly, (\ref{4}) and (\ref{5}) give
\BD
 m+m^*+k_2mnp=(1+k_3n)qp+m^*k_1np=qp+k_3qnp+m^*k_1np.
\ED
This, however, means
\BD
\frac{m+m^*}{np}\in \frac qn +\Z.
\ED
By (\ref{1}), $S(m,np)$ has the fractional part $q/n$.

{\em Case 2.} Let $n$ be even. Since $(q,n)=1$, $q$ is an odd number.
We choose a prime $p$ such that
\BE
\label{6}
   \CMOD 1{qp}n.
\EE
Note that it suffices to show one of the assertions $S(m,n')\in q/n+\Z$ or $S(m,n')\in -q/n+\Z$ for some $n'$, since
$S(-m,n')=-S(m,n')$ (see \cite[p. 26]{RaGr}). Accordingly, we may assume that $\CMOD q14$. If $\CMOD n04$,
condition (\ref{6}) implies $\CMOD p14$, in the case $\CMOD n24$, it only means that $p$ is odd, hence we
may choose $p$ such that $\CMOD p14$. Therefore, there exists an integer $m$ such that
\BE
\label{7}
  \CMOD m1{2n}\ \MB{ and }\  \CMOD{m^2}{-1}p.
\EE
From (\ref{6}) we obtain
\BD
  m^2+1\equiv 2\equiv 2qp \MB{ mod } 2n,
\ED
i. e., $m^2+1=2qp+2kn$, $k\in \Z$.
By (\ref{7}), $p$ divides $m^2+1$, so we have $m^2+1=2qp+2k_1np$, $k_1\in\Z$. Since $(m, 2np)=1$,
there exists an integer $m^*$ such that $mm^*=1+2k_2np$, $k_2\in\Z$. In addition, $m^*\equiv m\equiv 1$ mod $2n$.
Hence we write $m^*=1+2k_3n$, $n\in\Z$. Altogether, we obtain
\BD
 m^*(m^2+1)=m(1+2k_2np)+m^*=(1+2k_3n)(2qp+2k_1np),
\ED
which shows $m+m^*\in 2qp+ 2np\Z$, or
\BD
 \frac{m+m^*}{2np}\in \frac qn +\Z.
\ED
By (\ref{1}), $S(m, 2np)$ has the desired property.

\medskip
\noindent
{\em Example.} Let $n=132=3\cdot 4\cdot 11$ and $q=7$. Obviously, the fractional part of
$(m+m^*)/n$, $(m,n)=1$, lies in $(1/66)\cdot \Z$, hence it cannot be $7/132$.
Accordingly, we proceed as in the proof. Since $\CMOD q34$ we work with $-q=-7$
instead. Our prime $p$ has to satisfy $\CMOD 1{-7p}{132}$, i. e., $\CMOD p{113}132$. The smallest possible prime
of this kind is $p=509$. Now $m$ must be chosen such that
\BD
   \CMOD m1{264} \MB{ and }\CMOD{m^2}{-1}509.
\ED
A suitable number $m$ in the range $0\le m\le 2np-1$ is $m=133057$. Indeed, since $2np=134376$, we obtain
$S(m,2np)=-120-7/132$ and $S(-m,np)=S(n-m,np)=S(1319,134376)=120+7/132.$


\vspace{0.5cm}
\noindent
Kurt Girstmair            \\
Institut f\"ur Mathematik \\
Universit\"at Innsbruck   \\
Technikerstr. 13/7        \\
A-6020 Innsbruck, Austria \\
Kurt.Girstmair@uibk.ac.at

\end{document}